\documentclass[12pt]{article}
\usepackage{latexsym}
\usepackage{amssymb}
 
\newcommand{\proof}{{\noindent \bf Proof. }}
\newtheorem{thm}{Theorem}

\newtheorem{lem}{Lemma}

\newcommand{\x}{{\bf x}}

\newcommand{\C}{{\cal C}}

\newcommand{\D}{{\cal D}}
\newcommand{\F}{{\cal F}}
\newcommand{\G}{{\cal G}}

\newcommand{\Nat}{{\mathbb N}}

\date{}

\begin{document}
\begin{titlepage}
\title{\bf PATH SEPARATION BY SHORT CYCLES}
{\author{{\bf  G\'erard Cohen}
\\{\tt cohen@telecom-paristech.fr}
\\ENST
\\ FRANCE
\and{\bf Emanuela Fachini}
\\{\tt fachini@di.uniroma1.it}
\\''La Sapienza'' University of Rome
\\ ITALY
\and{\bf J\'anos K\"orner}
\thanks{Department of Computer Science, University of Rome, La Sapienza, 
via Salaria 113, 00198 Rome, ITALY}
\\{\tt korner@di.uniroma1.it}
\\''La Sapienza'' University of Rome
\\ ITALY}} 

\maketitle
\begin{abstract}

Two Hamilton paths in $K_n$ are separated by a cycle of length $k$ if their union contains such a cycle. For $k=4$ we bound the asymptotics 
of the maximum cardinality of a family of Hamilton paths in $K_n$ such that any pair of paths in the family is separated by a cycle of length $k.$
We also deal with related problems, including directed Hamilton paths.
\bigskip



\end{abstract}
\end{titlepage}

\section{Introduction}

The union of two graphs with the same vertex set is the graph on their common vertex set having as edge set the union of the edge sets of the two. The union of  two distinct Hamilton paths in $K_n$ always contains a cycle. We call a family of such paths an {\em odd family} if the union of any two of its members contains an odd cycle. 
It was observed in \cite{KMS} that the maximum cardinality of an odd family is only exponential in $n$, and precisely equal to the middle binomial coefficient  
${n \choose \lfloor n/2 \rfloor}$ if $n$ is odd while it is ${1 \over 2} {n \choose n/2}$ if $n$ is even. In the same way, we will call a family of Hamilton paths in $K_n$ an 
{\em even family} if the union of any pair of its members contains at least one even cycle. It is very easy to see that an even family can be considerably larger than an 
odd family and in fact can contain only an exponential factor times less than all Hamilton paths. This is obvious by realising that all the Hamilton 
paths contained in a fixed complete bipartite graph $K_{\lfloor n/2 \rfloor, \lceil n/2 \rceil}$ form an even family. Let $M(n, 2\Nat)$ be the largest cardinality of an even family of Hamilton paths from $K_n.$ The paper \cite{KMS} tightens the obvious bounds 
by showing that

\medskip
{\bf Theorem}\cite{KMS}

One has
$$\frac{n!}{2\Big [n-3+{\lceil\frac{n}{2}\rceil+1 \choose 2 }\Big ]} \leq M(n, 2\Nat)\leq \frac{n!}{2n}$$
if $n$ is odd, and
$$\frac{n!}{2{\frac{n}{2}+1 \choose 2 }} \leq M(n, 2\Nat)\leq \frac{n!}{8}$$
if $n$ is even.

\medskip
The problems just mentioned belong to the following general framework. Let $\F$ and $\D$ be two (not necessarily disjoint) families of graphs on the same vertex set 
$[n]$. We are interested in the largest cardinality $M(n, \F,\D)$ of a subfamily $\G \subseteq \F$ for which the union of any two 
different member graphs of $\G$ is in $\D$. A basic problem is to understand how the relationship of $\D$ to $\F$ influences the growth rate of  $M(n, \F,\D)$ as a function of 
$n.$ We do not have a real understanding of when it is that this growth rate is only exponential. The present paper is a continuation of the work of K\"orner,  Messuti and Simonyi \cite{KMS}. Throughout this paper, $\F$ will be the family of Hamilton paths in $K_n.$ This family is particularly interesting because of its close relationship with permutations of $[n].$ (As a matter of fact, a permutation of $[n]$ can be thought of as a (consecutively) oriented Hamilton path in the symmetrically complete directed graph on $n$ vertices).

Recently, the combinatorics of permutations became a growing and quite popular research topic. There is a good number of recent papers on intersection theorems, especially for permutations. The general topic of intersection theorems goes back to the seminal paper of Erd\H os, Ko and Rado \cite{EKR}. A famous and beautiful conjecture of Simonovits and S\'os \cite{Siso} was solved by Ellis, Filmus and Friedgut \cite{EFF}. Simonovits and S\'os \cite{Siso} have conjectured that the largest family of subgraphs of $K_n$ for which any two members of the family have a triangle in their intersection is obtained by considering all the subgraphs containing a fixed triangle and this maximum is unique. A remarkable feature of this problem, discovered 
in \cite{EFF} is that the family does not get larger if instead of a triangle, the pairwise intersection must only contain an odd cycle. 
This problem can be stated within our framework, in terms of the complements of the  graphs considered in \cite{Siso}, and replacing intersection by union. To do this, let 
$\F$ be the family of all graphs with vertex set $[n]$ and $\D$ be the family of those graphs on the same vertex set that have stability number (maximum cardinality of an independent set) at least 3.

In analogy to this, it was asked in \cite{KMS}
how large an odd family of Hamilton paths in $K_n$ can be if the union of any pair of its distinct members must contain a triangle (rather than just an arbitrary odd cycle). Let the largest 
cardinality of such a family be denoted by $M(n,3).$ Clearly, $M(n,3)$ is upper bounded by the maximum size of an odd family. But, can this bound be tight? It was shown in 
$\cite{KMS}$ that $M(5,3)=10$ verifying that at least for $n\leq 5$ the upper bound is tight. This gives, by an easy product construction, the lower bound 
$$M(n,3)\geq 10^{\lfloor n/5 \rfloor}.$$
No further progress has been made for odd cycles. Here we intend to ask the analogous question for even families.

\section{Two--part cycles}

We have recalled the fact that the largest cardinality of an even family of Hamilton paths in $K_n$ is less only by a quadratic factor than the total number $\frac{n!}{2}$ of 
Hamilton paths in the complete graph on $n$ vertices \cite{KMS}. In what follows we shall determine the asymptotics of $\hat{M}(n,4),$ the maximum cardinality of a family of Hamilton paths in $K_n$ such that their pairwise union contains a cycle of length four, and this cycle is the union of two subpaths from the two respective paths. As we shall see, $\hat{M}(n,4)$ grows, roughly speaking, only as the square root of $n!.$ More precisely,

\begin{thm}\label{thm:four}
We have
$$\lfloor n/2\rfloor!\leq \hat{M}(n,4) \leq (2\sqrt{3})^n\lceil n/2 \rceil!  \;.$$

\end{thm}

\proof

The lower bound construction is inspired by an idea from \cite{KM}. We will show that a balanced complete bipartite graph $K_{\lfloor n/2 \rfloor, \lceil n/2 \rceil}$ has 
at least this many Hamilton paths with the property that the pairwise union of any two of them contains a cycle of length four. To this end, let $A$ and $B$ be the two classes of the bipartition of  $K_{\lfloor n/2 \rfloor, \lceil n/2 \rceil}.$ Suppose that $|A|=\lfloor n/2 \rfloor.$ Let us fix an order of the elements of $B.$ To any permutation of the elements 
of $A$ we associate a Hamilton path in the bipartite graph as follows. The first vertex of each of these paths is the first element of the fixed order of $B$. The path alternates the 
elements of $A$ and $B$ where the various elements of $B$ appear in the fixed order, while the elements of $A$ appear in the order of the permutation of $A$ we consider.
Let now $\rho$ and $\sigma$ be two different permutations of $A.$ We claim that the union of the corresponding paths contains a cycle of length four. To see this, let $\pi(i),
i=1,\dots, n$ denote the image of $i\in [n]$ by the permutation $\pi.$ Then in the Hamilton path associated to $\rho,$ the vertex $\rho(i)$ is adjacent to the $i$'th and the 
$(i+1)$'st element of the fixed order of $B.$ Thus, for every $i\in [n].$ the vertices $\rho(i)$ and $\sigma(i)$ have the same pair of adjacent vertices in their respective Hamilton paths. Since $\rho$ and $\sigma$ are different, there is an $i\in [n]$ for which $\rho(i)\not =\sigma(i).$ Let us consider these two vertices along with their common pair of neighbors. These four vertices form a cycle of length four in $K_{\lfloor n/2 \rfloor, \lceil n/2 \rceil}.$

Let us turn to the upper bound. Let $\tau$ be an arbitrary Hamilton path from $K_n.$ To this path we associate a class of permutations as follows. Let 
$\tau$ be represented as the sequence $\tau(1),\tau(2), \dots, \tau(n).$ For the sake of simplicity suppose that $n$ is a multiple of 6 and let us partition the sequence of vertices of 
$\tau$ into consecutive disjoint groups of 6 elements, each. Let further each group be partitioned into disjoint consecutive pairs, called the {\em opening pair}, the {\em middle pair} and the 
{\em concluding pair}, respectively. We will say that two consecutive vertices, belonging to the same pair are {\em glued together}, since the partition into pairs is the same in all permutations from the same class. Let us associate with $\tau$ those permutations in which the opening pairs can be permuted arbitrarily among themselves, and the same is true 
for middle pairs and concluding pairs, too. More precisely, within all pairs, the order of the two elements of a pair remains the same, but the pairs of positions of opening pairs 
can be permuted in an arbitrary manner, and the same is true for middle pairs and concluding pairs, too. Let us denote the class of permutations associated with $\tau$ as 
$C(\tau).$ The sets $C(\cdot)$ of the various permutations represent a partition of the set of all the permutations of $[n].$ Each of the classes of this partition contains 
$[(n/6)!]^3$ permutations. We will refer to these sets as {\em filter set}s. We claim that the pairwise union of the Hamilton paths corresponding to two permutations from the same filter set does not contain any cycle of length four in which the edges from each path occur consecutively. 

If the union of two arbitrary paths contains a cycle of length four, this can happen in three different ways. In the first two cases the cycle is the union of two paths; one from each of the two Hamilton paths. 
In one case the two paths are of equal length, or, in the other case, it is the union of a path of one edge with a path of 3 edges. Finally, it can happen that the cycle alternates between edges from the two Hamilton paths. We have to show that neither of the first two cases occurs if the two paths are from the same filter set.

{\em Two paths of equal length}

Suppose that the union of the Hamilton paths $\tau$ and $\phi$ from the same filter set has a cycle of length four and this cycle is the union of two paths of two edges each 
from the two Hamilton paths $\tau$ and $\phi$. Therefore the two paths of two edges have the same endpoints, $\{a,b\}\in {[n] \choose 2}.$ This implies that the vertices $a$ and $b$ are both adjacent to some $x \in [n]$ in $\tau$ and to some $y \in [n]$ in $\phi,$ with $x\not=y.$ In particular, since the vertices of the paths come in pairs glued together in all the permutations, we see that $x$ is glued together with either $a$ or $b.$ The same is true for $y.$ Without loss of generality, suppose that in $\tau$ $x$ is glued together with $a$ and thus the same is true in $\phi$ as well, contradicting the hypothesis that $x$ and $y$ are different.

{\em Two paths of different length}

Suppose now that the union of the Hamilton paths $\tau$ and $\phi$ from the same filter set has a cycle of length four and this cycle is the union of one edge from $\tau$ with a 
path of three edges from $\phi.$ Let the common endpoints of the two paths be $a \in [n]$ and $b \in [n].$ Clearly, $a$ and $b$ cannot be glued together, otherwise they would 
be adjacent not only in $\tau$, but also in $\phi,$ a contradiction. Without loss of generality, we assume that $a$ precedes $b$ in $\tau.$ This implies that for some $i$ and $j$ we have $a=\tau(6i+j)$ where $2\leq j \leq 6$ and $j$ is even. Further, if 
$j=6$ then $b=\tau(6(i+1)+1)$, otherwise $b=\tau(6i+j+1).$ Note further that if $a=\tau(6i+j)$ then $a=\phi(6l+j)$ for some $l\in [n].$ 
Hence, $\phi^{-1}(a) \equiv \tau^{-1}(a) \; (\rm{mod} \;6).$
Likewise, we have 
$\phi^{-1}(b)\equiv \tau^{-1}(b) \; (\rm{mod}  \;6).$ Hence, if $j=6$, we have $a=\phi(6k)$ for some $k$ and $b=\phi(6m+1)$ for some $m\in [n],$ implying that $b$ cannot be 3 edges away from $a$ in $\phi.$ If $j<6$, we have $a=\phi(6l+j)$ and $b=\phi(6m+j+1)$ for some $m\in [n]$ and again, the two cannot be 3 edges apart in 
$\phi.$

Let us remark that we cannot exclude the presence (in the union of two paths from the same filter set) of alternating 4--cycles, i. e. cases in which the union of $\tau$ and $\phi$ has a 4-cycle whose edges are two vertex--disjoint edges from each of $\tau$ and $\phi.$ We have observed that all the filter sets have the same cardinality, 
$[(n/6)!]^3.$ In addition, the filter sets partition the family of all Hamilton paths from $K_n.$ Furthermore, a family of Hamilton paths no two of whose 
member graphs have two subpaths composing a 4--cycle cannot contain more than one path from any filter set. This gives
$$\hat{M}(n,4)\leq \frac{n!}{[(n/6)!]^3}=\frac{(n/2)!}{[( n/6 )!]^3}\cdot\frac{n!}{(n/2)!}\leq (2\sqrt{3})^n(n/2)!.$$

\hfill$\Box$

\section{Alternating cycles}

Let $A(n, 4)$ be the maximum cardinality of a family of Hamilton paths in the complete graph $K_n$ such that the union of any two of the paths from the family contains a 
cycle of length 4 with the restriction that the cycle contains two vertex--disjoint edges from both of the paths. We call such 4--cycles alternating. We do not know whether $A(n, 4)$ grows super-exponentially. 
The problem has an interesting connection with the question of a {\em reversing family} of permutations defined in \cite{El}. A family of permutations of $[n]$ is called 
reversing if for any pair of different permutations from the family there exists a pair $\{i,j\} \in {[n] \choose 2}$ of coordinates featuring the same two numbers but in different order. Let us denote by $R(n)$ the maximum cardinality of a reversing family of permutations of $[n].$ K\"orner conjectured (cf. the last section in the arXiv version of \cite{El}) that $R(n)$ grows only exponentially 
in $n.$ Obviously, $R(n)\geq 2^{\lfloor n/2 \rfloor}.$ This can be improved a bit, but the conjecture remains open. Cibulka \cite{C} has proved, improving an earlier bound of 
F\"uredi, Kantor, Monti and Sinaimeri \cite{FKMS} that 
$$ R(n) \leq  n^{n/2+O(\frac{n}{\log n})}.$$
(A family satisfying the same condition is called full of flips in \cite{FKMS}.) 
Our present interest in this question is motivated by the following observation.
\begin{lem}\label{lem:flip}

We have

$$A(n,4)\geq 2^{\lfloor {n-2 \over 4} \rfloor }R(\lfloor n/2 \rfloor).$$
\end{lem}

\proof
For the sake of simplicity suppose first that $n$ is even.
Let us fix a complete bipartite subgraph $K_{n/2, n/2}$ of $K_n.$
 Let the two disjoint stable sets of size $n/2$ of  $K_{n/2, n/2}$ be 
$A \subset [n]$ and $B\subset [n].$ Let us fix a linear order of the elements in $A$, and consider an arbitrary linear order of the elements in $B.$ 
To every linear order of $B$ we associate a perfect matching in $K_n$ the edges of which are connecting the $i$'th element of the fixed linear order of $A$ with the $i$'th element of the linear order of $B.$ Let us now consider a reversing family of permutations of $B$ of maximum cardinality 
$R(n/2).$ We observe that the pairwise unions of the corresponding perfect matchings contain an alternating 4--cycle. Let $\rho$ and $\sigma$ be two arbitrary permutations of $B$ from our reversing family. Let us suppose that a flip occurs in positions $\{a,b\}\in {B \choose 2}.$ Then the vertices of $B$  in these positions, alongside with the two vertices in $A$ in the corresponding positions define an alternating 4--cycle. 

To each of these perfect matchings we will associate $2^{\lfloor {n-2 \over 4} \rfloor }$ Hamilton paths in a way to satisfy our claim. Each of the Hamilton paths associated to a perfect matching will contain it. This guarantees already that two Hamilton paths associated with different perfect matchings have an alternating 4--cycle in their union. Let us now fix an arbitrary perfect matching from our family. 
Any orientation of the edges of a perfect matching  uniquely determines a directed Hamilton path in which the edges of the matching appear in the given order and with the given orientation. We group together the first and the second edge, and proceed in this way for every consecutive (disjoint) pair of edges. In any group of two edges we just defined we orient the first edge from $A$ to $B$.  For any sequence 
$\x \in \{-1,+1\}^{\lfloor n/4 \rfloor}$ we orient the second edge in the $i$'th pair of edges from $A$ to $B$ if 
$x_i=1$ and from $B$ to $A$ if  $x_i=-1.$ Actually, we only consider sequences $\x$ with their last coordinate fixed, say equal to 1.Then any pair of Hamiltonian paths associated with the same matching will generate an alternating 4-cycle having as vertices the 2 vertices of the second edge of 
the $j$'th couple together with the vertex from $B$ of the first edge in the couple and the vertex from $A$ of the first edge in the $(j+1)$'th couple, as soon as their defining sequences from $\{-1,+1\}^{\lfloor n/4 \rfloor}$ differ in the $j$'th coordinate. (Necessarily, the $j$'th coordinate is not the last one.)

For $n$ odd, suppose without loss of generality that $|A|=\lceil n/2 \rceil$ and fix the last vertex of the linear order of $A$ as the last vertex of each of the Hamilton paths.
Then repeat the same argument as before for the rest of the graph.

\hfill$\Box$

\section{Arbitrary 4--cycles}

The ideas developed in the previous two sections yield a proof of the following

\begin{thm}\label{thm:sug}

We have

$$\lfloor n/2\rfloor!\leq M(n,4) \leq n^{{3n \over 4}+O(\frac{n}{\log n})}$$

\end{thm}

\proof

The lower bound is immediate from the lower bound in Theorem \ref{thm:four}. In order to obtain the upper bound, we use the same set--up as in the proof of Theorem 
\ref{thm:four}. In particular, the proof is based on the concept of filter sets from said theorem.  

Consider a set of maximum cardinality among families of Hamilton paths satisfying the condition that their pairwise union contains a 4--cycle. Let us partition this set according to the filter sets these paths belong to.
We have shown in the proof of Theorem \ref{thm:four} that if the union of two Hamilton paths from the same filter set contains a 4--cycle, then this must be alternating 
between the two paths in the sense of the last section. Further, it is clear that no edge in such a 4--cycle contains two vertices that are glued together. The edges of these Hamilton paths define a bipartite graph with vertex set $[n].$ The two classes of the bipartition of the vertices are defined by the stable set formed by the $n/2$ first vertices of our glued--together pairs on the one hand, and the $n/2$ second vertices on the other. 

Notice that 
if we suppress the edges between glued vertices in any of our Hamilton paths, the rest forms a perfect matching between the same two classes of vertices of cardinality 
$n/2.$ Every such perfect matching in a fixed filter class fully determines the whole path. We regard these as permutations of the elements of a set of $n/2$ elements. 
Then an alternating 4--cycle in the union of two paths from the same filter set implies the existence of two vertices from both classes of the underlying partition (the very vertices of this cycle) such that they represent  a flip.
As observed  in the proof of Lemma \ref{lem:flip}, these permutations form a reversing family. Hence, by Cibulka's bound \cite{C} their number is at most $n^{n/4+O(\frac{n}{\log n})}.$ Considering that the number of filter sets is less than
${n \choose n/2}2^{n/2}({n \over 2})!<2^{2n}({n \over 2})!$, the upper bound follows.

\hfill$\Box$

\section{Directed graphs}

In the case of directed graphs new and interesting problems arise. A Hamilton path in a directed graph is a path in which the edges are oriented consecutively, meaning that 
every in--degree and out--degree of the path is at most 1. In other words, these can be considered as representations of permutations of $[n].$ The shortest cycle in a directed graph consists of two oppositely oriented arcs and therefore has length two. It is particularly interesting to examine the largest cardinality $C(n,2)$ of a family of Hamilton paths from the symmetrically complete digraph $\hat{K}_n$ such that the union of any two members of the family has a cycle of length two. It is somewhat surprising in comparison with our previous result that the growth rate of $C(n,2)$ is only exponential in $n.$ We have
\begin{thm}\label{thm:cy}

$$2^{\lfloor n/2 \rfloor} \leq C(n,2)\leq 2^{n-1}$$
\end{thm}

\proof

Let us start by proving the upper bound. We shall make use of a seminal result of Szele \cite{Sz} concerning Hamilton paths in tournaments. 
He proved the existence of a tournament $T_n$ on $n$ vertices which contains $\frac{n!}{2^{n-1}}$ well--oriented Hamilton paths. (We call a path well--oriented if all the in--degrees and out--degrees of its vertices are at most 1). Obviously, in this 
family of Hamilton paths there are no two whose union has a two--cycle, since in a tournament every edge has just one orientation. By symmetry, 
any Hamilton path in  $\hat{K}_n$ is contained in the same number of isomorphic copies of $T_n.$ Let $\C_n$ be a family of Hamilton paths satisfying our pairwise union condition. The foregoing observation implies that the density of $\C_n$ in the family of all the Hamilton paths from $\hat{K}_n$ cannot exceed the density of our family 
$\C_n$ within each of the copies of the tournament $T_n.$ This gives the upper bound.

For the lower bound let us consider the matching with edges ${2i-1, 2i}, i=1, \dots, {\lfloor n/2 \rfloor}$. Now, to this matching every 
sequence from $\x \in \{-1,+1\}^{\lfloor n/2 \rfloor}$ associates a unique directed Hamilton path in which the edges of the matching appear in their fixed order, as follows: the edge connecting $2i+1$ to $2i+2$ is going from the former to the latter if $x_i=+1$ and it goes in the reverse direction otherwise. The remaining edges of the corresponding directed Hamilton path are uniquely determined by the condition of making the whole path well--oriented, with the additional restriction that the edges of the matching should appear in their fixed order. It is clear that paths corresponding to different sequences from $\x \in \{-1,+1\}^{\lfloor n/2 \rfloor}$ have a two--cycle in their union. Actually, such two--cycles arise within the edges of the fixed matching.

\hfill$\Box$

The lower bound in this theorem is not asymptotically optimal. Angelo Monti has found, by computer search, a better construction for $n=6$ leading to a slight exponential improvement in the lower bound. (Clearly, every small construction can be used for asymptotic lower bounds through its Cartesian powers).

\section{Acknowledgement}

We gratefully acknowledge the ongoing interest of Angelo Monti in this project. Angelo suggested the study of the directed case and cycles of length two. 

Particular thanks are due to an anonymous referee who pointed out that the ideas in the first two sections give rise to Theorem \ref{thm:sug} and corrected several errors in our proofs.

\end{document}